\documentclass[11pt]{amsart}

\usepackage{color}

\usepackage{latexsym}
\usepackage{amsfonts, amssymb}

\newenvironment{thm}{\subsection{}{\textbf {Theorem.}}\em}{}
\newenvironment{prop}{\subsection{}{\textbf {Proposition.}}\em}{}
\newenvironment{cor}{\subsection{}{\textbf {Corollary.}}\em}{}
\newenvironment{lem}{\subsection{}{\textbf {Lemma.}}\em}{}

\newenvironment{pf}{\noindent{\textbf {Proof.}}}
{\hfill\eop\bigskip}

\newenvironment{defn}{\subsection{}{\textbf {Definition.}}\em}{\smallskip}
\newenvironment{eg}{\subsection{}{\textbf {Example.}}}{\smallskip}

\newenvironment{rem}{\subsection{}{\textbf {Remark.}}}{\smallskip}

\newcommand\fM{\ensuremath{\mathfrak M}}
\newcommand\fN{\ensuremath{\mathfrak N}}

\newcommand\fS{\ensuremath{\mathfrak S}}

\newcommand\fX{\ensuremath{\mathfrak X}}

\newcommand\cA{\ensuremath{\mathcal A}}
\newcommand\cB{\ensuremath{\mathcal B}}
\newcommand\cC{\ensuremath{\mathcal C}}

\newcommand\cE{\ensuremath{\mathcal E}}

\newcommand\cG{\ensuremath{\mathcal G}}
\newcommand\cH{\ensuremath{\mathcal H}}
\newcommand\cI{\ensuremath{\mathcal I}}

\newcommand\cK{\ensuremath{\mathcal K}}

\newcommand\cM{\ensuremath{\mathcal M}}

\newcommand\cP{\ensuremath{\mathcal P}}
\newcommand\cQ{\ensuremath{\mathcal Q}}
\newcommand\cR{\ensuremath{\mathcal R}}
\newcommand\cS{\ensuremath{\mathcal S}}
\newcommand\cT{\ensuremath{\mathcal T}}
\newcommand\cU{\ensuremath{\mathcal U}}
\newcommand\cV{\ensuremath{\mathcal V}}

\newcommand\bbC{\ensuremath{\mathbb C}}

\newcommand\bbM{\ensuremath{\mathbb M}}
\newcommand\bbN{\ensuremath{\mathbb N}}

\newcommand\bbR{\ensuremath{\mathbb R}}

\newcommand\hilb{\ensuremath{\mathcal H}}

\newcommand\ol{\ensuremath{\overline}}

\newcommand\eop{{{\hfil \ensuremath \Box}}}

\newcommand\norm{\ensuremath {\Vert}}

\newcommand\bofh{\ensuremath{\cB ( \cH)}}

\newcommand{\abs}[1]{\lvert#1\rvert}

\DeclareMathOperator{\rk}{rk}

\usepackage{latexsym}
\usepackage{amsfonts}

\begin{document}

\baselineskip 14pt


%
%


\title{On selfadjoint extensions of semigroups of partial isometries}

\thanks{${}^1$ Research supported in part by ARRS (Slovenia)}
\thanks{${}^2$ Research supported in part by NSERC (Canada)}

\thanks{{\ifcase\month\or Jan.\or Feb.\or March\or April\or May\or 
June\or
July\or Aug.\or Sept.\or Oct.\or Nov.\or Dec.\fi\space \number\day,
\number\year}}
\author
	[J. Bernik]{{Janez Bernik${}^1$}}
\address
	{Fakulteta za matematiko in fiziko\\
	Univerza v Ljubljani\\
	Jadranska 19 \\
	1000 Ljubljana \\
	Slovenia}
\email{janez.bernik@fmf.uni-lj.si}

\author
	[L.W. Marcoux]{{Laurent W.~Marcoux${}^2$}}
\author
	[A.I. Popov]{{Alexey I. Popov${}^2$}}
\author
	[H. Radjavi]{{Heydar Radjavi${}^2$}}
	
\address
	{Department of Pure Mathematics\\
	University of Waterloo\\
	Waterloo, Ontario \\
	Canada  \ \ \ N2L 3G1}
\email{LWMarcoux@uwaterloo.ca}
\email{a4popov@uwaterloo.ca}
\email{hradjavi@uwaterloo.ca}

\begin{abstract}
Let $\cS$ be a semigroup of partial isometries acting on a complex, infinite-dimensional, separable Hilbert space.  In this paper we seek criteria which will guarantee that the selfadjoint semigroup $\cT$ generated by $\cS$ consists of partial isometries as well.    Amongst other things, we show that this is the case when the set $\cQ(\cS)$ of final projections of elements of $\cS$ generates an abelian von Neumann algebra of uniform finite multiplicity.
\end{abstract}

\keywords{Partial isometry, semigroup, self-adjoint, abelian von Neumann algebra, multiplicity}
\subjclass[2010]{Primary: 47D03. Secondary: 47A65, 47B40, 20M20, 46L10}


\maketitle
\markboth{\textsc{  }}{\textsc{}}

\section{Introduction}
A \textbf{semigroup} is a non-empty set $\cS$ equipped with an associative binary operation $\circ: \cS \times \cS \to \cS$.   Semigroups are   pervasive  in mathematics, and they appear in a multitude of contexts, algebraic, geometric and analytic.  After groups, perhaps the most important class of semigroups is that of \textit{inverse semigroups}.   Inverse semigroups were originally discovered by Wagner (sometimes written ``Vagner" in the literature) in 1952.   He referred to them as ``generalized groups".   They were independently discovered by Preston in 1954, and it is to him that we owe their present name.  

An element $s$ of a semigroup $\cS$ is said to be \textbf{regular} if there exists $t \in \cS$ such that $s t s = s$.  The semigroup $\cS$ is said to be \textbf{regular} if each of its elements is regular.   Note that when $s \in \cS$ is regular, we have that $s (tst) s = (s t s) t s = s t s = s$, while $(t s t) s (t s t) = t ((s t s) t s) t   = t ((s) t s) t = t s t$, so that $tst$ and $s$ are (said to be) \textbf{inverses} of each other.   (In other words, regular elements admit inverses.)  \textbf{Inverse semigroups} may be defined as semigroups in which every element admits a \textit{unique} inverse, or equivalently (see~\cite{Pet1984}), inverse semigroups are regular semigroups in which the idempotent elements commute.

In the same way that Cayley's Theorem allows one to view each group $\cG$ as a group of permutations acting on a set (which we may take to be $\cG$ itself), the Wagner-Preston Theorem~\cite{Pet1984, Pre1954} identifies inverse semigroups as particular subsemigroups of the \textbf{symmetric inverse semigroup} $\cI_X$ of all partial one-to-one transformations of a fixed set $X$.   (A \textbf{partial one-to-one transformation} on $X$ is a bijective map $\alpha: Y \to Z$, where $Y$ and $Z$ are subsets of $X$.)

In operator theory and operator algebras, semigroups have also been of recent interest.   The $C^*$-algebras studied by Cuntz~\cite{Cun1977}, by Cuntz and Krieger~\cite{CK1980}, the graph $C^*$-algebras studied by Kumjian, Pask, Raeburn (and many others) ~\cite{Rae2005},  the non-selfadjoint free semigroup algebras studied by Davidson, Katsoulis, Pitts, \emph{et al.}~\cite{DKP2001}, and the non-commutative disc algebras studied by Popescu~\cite{Pop1996} are just a few of the many places where semigroups have played a central role.  In each of the above cases, the semigroups have consisted of semigroups of partial isometries, which we now define.

Let $\hilb$ be a complex Hilbert space, which we shall normally assume to be infinite-dimensional and separable.   Denote by $\bofh$ the $C^*$-algebra of bounded linear operators acting on $\hilb$.    For any subset $\cA \subseteq \bofh$, let $\cA^* := \{ A^*: A \in \cA\}$.   A \textbf{partial isometry} $V \in \bofh$ is a map for which the operator $Q_V := V V^*$ is a selfadjoint idempotent (i.e., a \emph{projection}, called the \textbf{range projection} of $V$), or equivalently, for which  $P_V := V^* V$ is a projection (called the \textbf{initial projection} of $V$).    Given a collection $\cA$ of partial isometries,  let us also write $\cP(\cA) = \{ P_S : S \in \cA\}$ and $\cQ(\cA) = \{ Q_S: S \in \cA\}$.   If $\cS$ is a semigroup of partial isometries, then clearly $\cS^*$ is also a semigroup of partial isometries, and $\cP(\cS) = \cQ(\cS^*)$, while $\cQ(\cS) = \cP(\cS^*)$.  Given a commuting family of projections $\cP$, we denote by $BA(\cP)$ the Boolean algebra generated by these projections.

The notion of a partial isometry may profitably be thought of as a linearization of the concept of a partial one-to-one map, as a partial isometry $V$ acts as an isometric linear map (and in particular as a bijection) from $P_V \hilb$ to $Q_V \hilb$, and has kernel $(P_V \hilb)^\perp$.   It is elementary to check that $V^* V V^* = V^*$ and that $V V^* V = V$ for every partial isometry.   In general, however, the product of two partial isometries need not be a partial isometry.   Indeed, we have the following result of Halmos and Wallen to which we shall often have recourse.


\begin{thm} \emph{(}\textbf{Halmos-Wallen\ \cite{HW1969}}\emph{)}  \label{HalmosWallen}
Let $\hilb$ be a Hilbert space and suppose that $V$ and $W$ are partial isometries acting on $\hilb$.  The product $V W$ is a partial isometry if and only if $Q_W$ commutes with $P_V$.
\end{thm} 


\bigskip

Let $\cS$ be a selfadjoint semigroup  of partial isometries in $\bofh$.   Then $\cS$ is clearly  regular (with $V^*$ acting as an inverse to $V$ as seen above).   Furthermore, if $E \in \cS$ is an idempotent, then by the Halmos-Wallen Theorem above, $Q_E$ commutes with $P_E$, from which it is easily seen that $E=E^*$ and so $E$ is an orthogonal projection.   Finally, if $E$ and $F$ are idempotents in $\cS$, then $E=E^*$ and $F=F^*$ and therefore,  as $E F \in \cS$ is a partial isometry, a second application of the Halmos-Wallen Theorem shows us that $Q_F = F$ and $P_E = E$ must commute.   Thus $\cS$ is necessarily an inverse semigroup.   A result of Duncan and Paterson~\cite{DP1985} shows that the converse also holds.   Given an inverse semigroup $\cS$, and denoting the unique inverse of $t \in \cS$ by $t^*$, we let $\cS$ act upon $\ell_2(\cS)$ via Barnes' \textit{left regular representation}~\cite{Bar1976} 
\[
\begin{array}{rccc}
	\pi: & \cS & \to & \cB(\ell_2(\cS)) \\
	     & s & \mapsto & \pi(s),
\end{array} \]
where $\pi(s) \big( \sum_{t \in \cS} a_t t \big) = \sum_{t t^* \le s^* s} a_t s t$.  (Note that $t t^*$ and $s^* s$ are commuting idempotents in $\cS$.   We say that $t t^* \le s^* s$ if $(t t^*) (s^* s) = t t^*$.)  The Duncan-Paterson result shows that the resulting set $\pi (\cS) = \{ \pi(s) : s \in \cS \}$ is an inverse semigroup of partial isometries which is isomorphic (as an inverse semigroup) to $\cS$.

If  $\cA \subseteq \bofh$, we denote by $\cA^\prime$ the \textbf{commutant} of $\cA$; that is,  $\cA^\prime := \{ T \in \bofh: T A =  A T \mbox{ for all } A \in \cA \}$.   As a consequence of the Halmos-Wallen result above, we obtain the following.


\begin{rem} \label{rem00.5}
Suppose that $\cS$ is a semigroup of partial isometries.   Then $\cP(\cS) \subseteq \cQ(\cS)^\prime$ and $\cQ(\cS) \subseteq \cP(\cS)^\prime$.
\end{rem} 

\bigskip

The work below is a natural extension of the results obtained in~\cite{PR2013} by the last two authors.   Given a set $\cA$ of operators acting on $\hilb$, let $\langle \cA \rangle$ denote the semigroup generated by $\cA$.   The basic question that we wish to address here is, given a  semigroup $\cS$ of partial isometries, when is the selfadjoint semigroup $\cT = \langle \cS \cup \cS^* \rangle$ generated by $\cS$ a semigroup of partial isometries?

Of course, if a semigroup $\cS$ is both a maximal semigroup of partial isometries and $\cT = \langle \cS \cup \cS^* \rangle$ consists of partial isometries, then $\cS = \cT$, so that $\cS$ is already selfadjoint.  However, the set $\cS$ of all linear isometries on $\ell_2(\bbN)$, together with the zero operator, serves as an example of a maximal semigroup of partial isometries for which the corresponding selfadjoint semigroup $\cT$ does not consist exclusively of partial isometries.   This is not too surprising.   Indeed, it is clear from the comments above that if $\cT$ is to be a selfadjoint semigroup of partial isometries, then it must be an inverse semigroup, and hence all idempotents in $\cT$ must commute.  In particular, since $\cP(\cS)$ and $\cQ(\cS)$ consist of idempotents in $\cT$, at the very least we require that $\cP(\cS) \cup \cQ(\cS)$ be a commutative set of projections.  

Recall that a set $\cA \subseteq \bofh$ is said to act \textbf{irreducibly} on $\hilb$ (i.e. $\cA$ is \textbf{irreducible}) if the only closed, invariant subspaces for $\cA$ are the trivial spaces $\{0\}$ and $\hilb$.   Otherwise, we say that $\cA$ is \textbf{reducible}.   In terms of irreducible semigroups, we have the following result.  Although it is surely known, we have not been able to locate an explicit reference in the literature, and so we include a proof of it.


\begin{thm} \label{irreducible}
Let $\hilb$ be a Hilbert space and $\cS \subseteq \bofh$ be a semigroup. Then $\cS$ is irreducible if and only if
$$
A\cS B\ne\{0\}
$$
for every nonzero $A, B \in \bofh$.
\end{thm} 

\begin{pf}
It is clear that there is no loss of generality in assuming that the identity operator $I \in \cS$.

First suppose that there exist nonzero $A, B \in \bofh$ so that $A S B = 0$ for all $ S \in \cS$.   Choose $0 \not = x \in \hilb$ such that $B x \not = 0$ and $0 \not = y \in \hilb$ so that $A^* y \not = 0$.  Consider $\cM := \mathrm{span}\, \{ S B x : S \in \cS \}$. Since $\cS$ is a semigroup, $\cM$ is $\cS$-invariant subspace. Clearly $\cM \not = \{ 0\}$, because $B x \in \cM$.   If $m \in \cM$, then there exist $n \ge 1$, scalars $\alpha_1, \alpha_2, ..., \alpha_n \in \bbC$ and $S_1, S_2, ..., S_n \in \cS$ so that $m = \sum_{i=1}^n \alpha_i S_i B x$.   Now 
\begin{align*}
\langle m, A^* y \rangle
	&= \sum_{i=1}^n \alpha_i \langle S_i B x, A^* y \rangle \\
	&= \sum_{i=1}^n \alpha_i \langle A S_i B x, y \rangle \\
	&= \sum_{i=1}^n \alpha_i \langle 0 x, y \rangle = 0.
\end{align*}
Thus $\cM$ (and therefore $\ol{\cM}$) is contained in $(A^* y)^\perp$, so that $\ol{\cM}$ is a proper closed subspace of $\hilb$.   Hence $\cS$ is reducible.   

Conversely, suppose that $\cS$ is reducible, and let $\cM$ be a non-trivial, closed, invariant subspace for $\cS$.   Let $0 \ne x \in \cM$ and $0 \not = y \in \cM^\perp$ be vectors of norm one.   Let $B$ (resp. $A$) denote the orthogonal projection of $\hilb$ onto $\bbC x$ (resp. of $\hilb$ onto $\bbC y$).   If $S \in \cS$ and $z \in \hilb$, then 
\[
A S B z = \langle z, x \rangle A S B x =  \langle z, x \rangle A S x	= \langle z, x \rangle \ \langle S x, y \rangle y = 0 \]
as $S x \in \cM$ and $y \in \cM^\perp$.  Thus $A S B = 0$.   Since $S \in \cS$ was arbitrary, $A \cS B = \{ 0 \}$.
\end{pf}


Corollary~3.12 of~\cite{PR2013} states that if $\cS \subseteq \bofh$ is a semigroup of partial isometries which acts irreducibly on $\hilb$, and if $\cS$ contains a non-zero compact operator, then the selfadjoint semigroup $\cT$ generated by $\cS$ consists of partial isometries.   That this  result does not generally hold in the absence of the irreducibility hypothesis is demonstrated by the following example.


\begin{eg} \label{eg04.5}
Let $E = \begin{bmatrix} 1/2 & 1/2 \\ 1/2 & 1/2 \end{bmatrix}$ and $F = \begin{bmatrix} 1 & 0 \\ 0 & 0 \end{bmatrix}$, so that $E$ and $F$ are two non-commuting projections in $\bbM_2(\bbC)$.   Consider the partial isometries $A, B,$ and $C \in \bbM_8(\bbC)$ defined as follows: 
\begin{align*}
A &= \begin{bmatrix} 0 & E & 0 & 0  \\ 0 & 0 & 0 & 0 \\ 0 & 0 & 0 & 0 \\ 0 & 0 & 0 & 0  \end{bmatrix},  \\
B &= \begin{bmatrix} 0 & 0 & 0 & I_2  \\ 0 & 0 & 0 & 0 \\ 0 & 0 & 0 & 0 \\ 0 & 0 & 0 & 0  \end{bmatrix},  \\
C &= \begin{bmatrix} 0 & 0 & 0 & 0  \\ 0 & 0 & 0 & 0 \\ 0 & 0 & 0 & F \\ 0 & 0 & 0 & 0  \end{bmatrix}. 
\end{align*}
Observe that the product of any two elements of $\cS_0 := \{ A, B, C, 0\}$ is zero, from which we deduce that $\cS_0$ is a semigroup of partial isometries in $\bbM_8(\bbC)$.    We extend $\cS_0$ to a unital semigroup $\cS =  \cS_0 \cup \{ I_8\} $, where $I_8$ denotes the identity operator in $\bbM_8(\bbC)$.

Now $P_A = A^* A = \mathrm{Diag}(0, E, 0, 0)$,   $Q_A = A A^* = \mathrm{Diag}(E, 0, 0, 0)$, $P_B = \mathrm{Diag}(0, 0, 0, I_2)$, $Q_B = \mathrm{Diag}(I_2, 0, 0, 0)$,  $P_C = \mathrm{Diag}(0, 0, 0, F)$, $Q_C = \mathrm{Diag}(0, 0, F, 0)$, $P_{I_8} = Q_{I_8} = I_8$, so  that $\cP(\cS) \cup \cQ(\cS)$ is commutative. 
We claim that there is no selfadjoint semigroup $\cT$ of  $\bbM_8(\bbC)$ which contains $\cS$ and consists of partial isometries.   For suppose otherwise.    Since $\cT \supseteq \cS$ is  selfadjoint, it follows that  $\cP(\cS) \cup \cQ(\cS) \subseteq \cT$ and $B P_C B^* = \mathrm{Diag}(F, 0, 0, 0) \in \cT$ as well.  A simple calculation now shows that  $Q_A  B P_C B^* = \mathrm{Diag} (E F, 0, 0, 0) \in \cT$ is not a partial isometry, a contradiction.   Of course, we could also use an application of Zorn's Lemma (or a simple dimension argument) to extend $\cS$ to a maximal semigroup $\cS_{\textrm{max}}$ of partial isometries on $\bbC^8$, for which $\cT_1 := \langle \cS_{\textrm{max}} \cup \cS^*_{\textrm{max}} \rangle \not = \cS_{\textrm{max}}$, for otherwise $\cT$ above would consist of partial isometries.
\end{eg}


A relatively straightforward computation shows that the von Neumann algebra $W^*(\cQ(\cS))$ generated by $\cS$ is unitarily equivalent to $\mathrm{Diag} (\alpha, \beta, \gamma, \delta I_5)$, where $\alpha, \beta, \gamma, \delta \in \bbC$ and where $I_5$ denotes a $5 \times 5$ identity matrix.  In particular, it is clear that the multiplicity of $W^*(\cQ(\cS))$ as a von Neumann algebra is \emph{not uniform}.  We will show that if the von Neumann algebra $W^*(\cQ(\cS))$ is uniform and finite, then the semigroup $\cS$ can always be extended to a selfadjoint semigroup (compare with Theorem~\ref{thm2.24}).

It is clear that the semigroup $\cS$ from Example~\ref{eg04.5} does not act irreducibly on $\bbC^8$.   In Theorem~\ref{uniform-multiplicity}, we prove that if $\cS$ is an \emph{irreducible} (unital) semigroup  for which $\cQ(\cS)$ is commutative, then the multiplicity of $W^*(\cQ(\cS))$ is automatically uniform.     Still, as we shall see later, irreducibility is not the only issue.   We adapt C.~Read's construction of a semigroup of partial isometries~\cite{Rea2005}, initially developed to demonstrate that $\bofh$ is a free semigroup algebra,  to provide an example of an irreducible semigroup $\cS$ of partial isometries for which the selfadjoint semigroup $\cT$ generated by $\cS$  does not consist of partial isometries (c.f. Theorem~\ref{not-extendable}).     The other key obstacle to ``${}^*$-extendability" is whether the multiplicity of the von Neumann algebra is finite or infinite.

We also establish other conditions which guarantee that an irreducible semigroup $\cS$ of partial isometries is ${}^*$-extendable.  In particular, we shall see that this is the case when $\cP(\cS) \cup\cQ(\cS)$ admits a minimal non-zero projection (Theorem~\ref{thm08}).

We finish this section by pointing out that while a number of results below (in particular, those at the end of the third section) are stated in terms of conditions upon the set $\cQ(\cS)$ of range projections associated to the semigroup $\cS$,   similar result can be obtained by replacing $\cQ(\cS)$ by $\cP(\cS)$ and considering $\cS^*$ instead of $\cS$.
 

\section{General results}

Throughout the remainder of the paper, \emph{we shall implicitly assume that the semigroups $\cS$ under consideration are unital}.

\begin{lem} \label{lem01}
Let $\cS$ be a semigroup of partial isometries, and let $m \ge 1$.  
If $R_1$, $R_2$, $\dots$, $R_m$, $S \in \cS$, then 
\[
S  (Q_{R_1} Q_{R_2} \cdots Q_{R_m}) S^* = Q_{S R_1} Q_{S R_2} \cdots Q_{S R_m}. \]
\end{lem}

\begin{pf}
First observe that $S   = S P_S = S S^* S$, as $S$ is a partial isometry.  It follows from Remark~\ref{rem00.5} that $P_S$ commutes with each $Q_{R_k}$, $1 \le k \le m$, and since $P_S = P_S^j$ for all $j \ge 1$, we get
\begin{align*}
S  (Q_{R_1} Q_{R_2} \cdots Q_{R_m}) S^*
	&= S P_S^{m-1} (Q_{R_1}\  Q_{R_2}\  \cdots \  Q_{R_m}) S^* \\
	&= S Q_{R_1} P_S Q_{R_2} P_S \cdots P_S Q_{R_m} S^* \\
	&= S Q_{R_1} S^* S Q_{R_2} S^* S \cdots S Q_{R_m} S^* \\
	&= (S R_1 R_1^* S^*) \ (S R_2 R_2^* S^*) \ \cdots (S R_m R_m^* S^*) \\
	&= Q_{S R_1} \ Q_{S R_2} \ \cdots Q_{S R_m}. 
\end{align*}	 
\end{pf}

We shall also require the following variant of Lemma~\ref{lem01}.

\begin{lem} \label{lem01alternate}
Let $\cS$ be a semigroup of partial isometries and $R_1, R_2, ..., R_m \in \cS$.   Let $V \in \bofh$ be a partial isometry and suppose that $P_V \in \cQ(\cS)^\prime$.   Then 
\begin{enumerate}
	\item[(a)]
	$V (Q_{R_1} Q_{R_2} \cdots Q_{R_m}) V^* = Q_{V R_1} Q_{V R_2} \cdots Q_{V R_m}$, and 
	\item[(b)]
	If $\cQ(\cS)$ is abelian, then $V (Q_{R_1} Q_{R_2} \cdots Q_{R_m}) V^*$ is a projection.
\end{enumerate}
\end{lem}

\begin{pf}
\begin{enumerate}
	\item[(a)]
	The proof of this is identical to that of Lemma~\ref{lem01}.
	\item[(b)]
	Let $X = 	 V (Q_{R_1} Q_{R_2} \cdots Q_{R_m}) V^*$.  Since each term in the product which defines $X$ has norm at most one, $\norm X \norm \le 1$.   On the other hand, 
	\begin{align*}
	X^2 
		&=  V (Q_{R_1} Q_{R_2} \cdots Q_{R_m}) V^* \  V (Q_{R_1} Q_{R_2} \cdots Q_{R_m}) V^* \\
		&=  V (Q_{R_1} Q_{R_2} \cdots Q_{R_m}) P_V (Q_{R_1} Q_{R_2} \cdots Q_{R_m}) V^* \\
		&=  V P_V (Q_{R_1}^2 Q_{R_2}^2 \cdots Q_{R_m}^2) V^* \\
		&=  V (Q_{R_1} Q_{R_2} \cdots Q_{R_m}) V^* \\
		&=  X. 
	\end{align*}
	Thus $X$ is an idempotent of norm at most one, and hence a projection.
\end{enumerate}
\end{pf}


\begin{prop} \label{prop03}
Let $\cS$ be a semigroup of partial isometries, and suppose that  $\cQ(\cS)$ is commutative.  Let $T \in \cS$.   Then the semigroup $\cS_0$ generated by $Q_T$ and by $\cS$ (i.e. the set of all words in $Q_T$ and elements of $\cS$) consists of partial isometries, $BA(\cQ(\cS_0)) = BA (\cQ(\cS))$ and   $\cQ(\cS_0)$ is again a commuting family of projections.   
\end{prop}

\begin{pf}
Recall that we have assumed that $\cS$ is unital.  In this case, any word $W$ in $Q_T$ and $\cS$ may be expressed in the form:  
\[
W := S_m \ Q_T \ S_{m-1} \ Q_T \ \cdots \ Q_T \ S_1, \]
for some $m \ge 2$ and elements $S_1, S_2, ..., S_m \in \cS$. We shall show that 
\[
Q_W = Q_{S_m T} \ Q_{S_m S_{m-1} T}\  \cdots\  Q_{S_m S_{m-1} \cdots S_2 T}\  Q_{S_m S_{m-1} \cdots S_1}. \]
Since $\cQ(\cS)$ is commutative, this forces $Q_W$ to be a projection, from which it follows that $W$ is indeed a partial isometry, as required.   Furthermore, it shows that $Q_W \in BA (\cQ(\cS))$, so that $BA(\cQ(\cS_0)) \subseteq BA(\cQ(\cS))$.   Since $\cS \subseteq \cS_0$, the reverse containment is trivial, and so equality of these two Boolean algebras ensues.

We shall use a proof by induction on $m$.  We begin with $m=2$. 
If $W = S_2 Q_T S_1$ then by Lemma~\ref{lem01} we find that
	\begin{align*}
	Q_W 
		&= S_2 Q_T S_1  S_1^* Q_T S_2^* \\
		&= S_2 Q_T Q_{S_1}  Q_T S_2^* \\
		&= S_2 Q_T Q_{S_1}  S_2^* \\
		&= Q_{S_2 T} \ Q_{S_2 S_1 }.
	\end{align*}
Under the induction assumption that the statement holds for $m=n$, we shall prove it for $m = n+1$.
	
Suppose that $V = 	S_n \ Q_T \ S_{n-1} \ Q_T \ \cdots \ Q_T \ S_1$ implies that 
	\[
	Q_V = Q_{S_n T} \ Q_{S_n S_{n-1} T}\  \cdots\  Q_{S_n S_{n-1} \cdots S_2 T}\  Q_{S_n S_{n-1} \cdots S_1}. \]
	
	Let $W = S_{n+1} Q_T V$.   Then 
	\begin{align*}
	Q_W 
		&= S_{n+1} Q_T Q_V Q_T S_{n+1}^* \\
		&= S_{n+1} Q_T (Q_{S_n T} \ Q_{S_n S_{n-1} T}\  \cdots\  Q_{S_n S_{n-1} \cdots S_2 T}\  Q_{S_n S_{n-1} \cdots S_1}) S_{n+1}^* \\
		&= Q_{S_{n+1} T}\  Q_{S_{n+1} S_n T} \ \cdots \ Q_{S_{n+1} S_{n} \cdots S_1}.
	\end{align*}
	By induction, the result holds for all $m \ge 2$.

\end{pf}

\begin{rem} \label{rems03.5}

	By applying Remark~\ref{rem00.5} to the semigroup $\cS_0$ above, we see that $\cP(\cS_0) \subseteq \cQ(\cS_0)^\prime$ and $\cQ(\cS_0) \subseteq \cP(\cS_0)^\prime$.

\end{rem}


\begin{prop} \label{prop04}
Let $\cS \subseteq \bofh$ be a semigroup of partial isometries for which $\cQ(\cS)$ is commutative.   Then there exists a  semigroup $\cS_{\textrm{max}}$ of partial isometries which is maximal with respect to the conditions that 
\begin{enumerate}
	\item[(i)]
	$\cS_\textrm{max} \supseteq \cS$,  and 
	\item[(ii)]
	$W^*(\cQ(\cS_\textrm{max})) = W^*(\cQ(\cS))$.   
\end{enumerate}
Furthermore, $\cQ(\cS_{\textrm{max}}) \subseteq \cS_{\textrm{max}}$.
\end{prop}

\begin{pf}
Observe that condition (ii) implies that $\cQ(\cS_{\textrm{max}})$ is commutative.  Not surprisingly, we shall appeal to Zorn's Lemma.

Let 
\begin{align*}
\fS &:= \{ \cR \subseteq \bofh: \cR \supseteq \cS \mbox{ is a semigroup of partial isometries}, \\
	&\ \ \ \ \ \ \ \ \ \ \ \ \ \ \ \ \ \ \ \ 
	  \ \ \ \ \ \ \ \ \ \ \ \  \mbox{ and } W^*(\cQ(\cR)) = W^*(\cQ(\cS)) \}. 
\end{align*}	
Clearly $\fS$ can be partially ordered by inclusion, so that $\cR_1 \le \cR_2$ provided that $\cR_1 \subseteq \cR_2$.   Since $\cS \in \fS$, we see that $\fS \not = \varnothing$.   

If $\cC = \{ \cR_\beta: \beta \in \Gamma \}$ is a chain in $\fS$, then $\cR_0 := \cup_{\beta \in \Gamma} \cR_\beta$ is easily seen to be a semigroup of partial isometries containing $\cS$.   Furthermore, if $Q \in \cQ(\cR_0)$, then $Q = Q_T$ for some $T \in \cR_0$, and so $T \in \cR_{\beta_0}$ for some $\beta_0 \in \Gamma$.   Hence $Q = Q_T \in \cQ(\cR_{\beta_0})$.   But then $Q \in W^*(\cQ(\cR_{\beta_0})) = W^*(\cQ(\cS))$, so $W^*(\cQ(\cR_0)) \subseteq W^*(\cQ(\cS))$.   The reverse inclusion is trivial, given that $\cS \subseteq \cR_0$.    This shows that $\cR_0 \in \fX$ is an upper bound for $\cC$.    We may therefore apply Zorn's Lemma to conclude that $\fS$ admits a maximal element, say $\cS_{\textrm{max}}$.  

Suppose that there exists $T \in \cS_{\textrm{max}}$ such that $Q_T \not \in \cS_{\textrm{max}}$.   By applying Proposition~\ref{prop03} to the pair $T$ and $\cS_{\textrm{max}}$, we obtain a  semigroup $\cS_0 = \langle \cS_{\textrm{max}} \cup \{Q_T\} \rangle$ of partial isometries strictly containing $\cS_{\textrm{max}}$ for which $BA(\cQ(\cS_0)) = BA(\cQ(\cS_{\textrm{max}}))$, from which it follows that 
\[
W^*(\cQ(\cS_0)) = W^*(\cQ(\cS_{\textrm{max}})) = W^*(\cS), \]
contradicting the maximality of $\cS_{\textrm{max}}$.   Thus $\cQ(\cS_{\textrm{max}}) \subseteq \cS_{\textrm{max}}$.
\end{pf}

It is interesting to note that although we can extend the semigroup $\cS$ to a larger semigroup $\cS_{\textrm{max}}$ of partial isometries which includes all of its range projections, there is no reason to believe that $\cP(\cS_{\textrm{max}})$ is commutative, though by Remark~\ref{rem00.5} it is clear that 
\[
\cP(\cS_{\textrm{max}}) \subseteq \cQ(\cS_{\textrm{max}})^\prime \subseteq W^*(\cQ(\cS_{\textrm{max}}))^\prime  = W^*(\cQ(\cS))^\prime.\]
As a consequence of this, if $W^*(\cQ(\cS))$ happens to be a maximal abelian selfadjoint algebra (i.e. a masa) in $\bofh$, then $W^*(\cQ(\cS)) = W^*(\cQ(\cS))^\prime$ and $\cP(\cS_{\textrm{max}}) \cup \cQ(\cS_{\textrm{max}})$ will be a commuting set of partial isometries.     Suppose that $T \in \cS_{\textrm{max}}$.  Then by Proposition~\ref{prop03}, the semigroup $\cR$ generated by $Q_{T^*}$ and $\cS_{\textrm{max}}^*$ consists of partial isometries and $W^*(\cQ(\cR)) = W^*(\cQ(\cS_{\textrm{max}}^*)) = W^*(\cP(\cS_{\textrm{max}}))$.   Upon taking adjoints, we find that $\cR^*$, which is the semigroup generated by $P_T$ and $\cS_{\textrm{max}}$, consists of partial isometries.  But $\cQ(\cS_{\textrm{max}}) \subseteq \cS_{\textrm{max}}$ implies that $\cQ(\cS_{\textrm{max}}) \subseteq \cP(\cS_{\textrm{max}})$.      Thus the facts that $\cQ(\cR^*) \subseteq \cP(\cR^*)^\prime$ and $\cS_{\textrm{max}} \subseteq \cR^*$ imply that 
\[
\cQ(\cR^*) \subseteq \cP(\cS_{\textrm{max}})^\prime \subseteq \cQ(\cS_\textrm{max})^\prime = W^*(\cQ(\cS_\textrm{max}))^\prime = W^*(\cQ(\cS))^\prime = W^*(\cQ(\cS)).\] 
From this it easily follows that $W^*(\cQ(\cR^*)) = W^*(\cQ(\cS))$.   Since $\cS_{\textrm{max}}$ is a maximal semigroup of partial isometries satisfying the conditions of Proposition~\ref{prop04}, it follows that  $\cR^* = \cS_{\textrm{max}}$ and hence that $P_T \in \cS_{\textrm{max}}$.   Thus $\cP(\cS_{\textrm{max}}) \subseteq \cS_{\textrm{max}}$.    Hence $\cP(\cS_{\textrm{max}}) \subseteq \cQ(\cS_{\textrm{max}})$ (as each projection in a semigroup is both its own initial and final projection), and so $\cP(\cS_{\textrm{max}}) = \cQ(\cS_{\textrm{max}})$ when $W^*(\cQ(\cS))$ generates a masa.

Proposition~\ref{max-prop} of the next section should be seen as a generalization of this idea.


\begin{thm} \label{thm05}
Suppose that $\cS$ is a semigroup of partial isometries and that $\cP(\cS) = \cQ(\cS)$.   Then the semigroup $\cT$ generated by $\cS$ and $\cS^*$ consists of partial isometries.
\end{thm}

\begin{pf}
Again, we recall that we are assuming that $\cS$ is unital.

Suppose that $R \in \cP(\cS) = \cQ(\cS)$, and write $R = S^* S = T T^*$, where $S, T \in \cS$.  Observe that for any $V \in \cS$, $VT, S V  \in \cS$ and so 
\begin{align*}
V R V^* &= V T T^* V^* = Q_{V T} \in \cQ(\cS) \mbox{ and } \\
V^* R V &= V^* S^* S V = P_{S V} \in \cP(\cS). 
\end{align*}

Furthermore, since $I \in \cS$, any element $W$ in the semigroup $\cT$ may be expressed as 
\[
W = S_1 T_1^* S_2 T_2^* ... S_n T_n^* \]
for some choice of $S_1, S_2, ..., S_n, T_1, T_2, ..., T_n \in \cS$.    We will show that $W$ is a partial isometry, which will complete the proof.   

More precisely, we shall use induction (on $n$) to show that $WW^* \in \cQ(\cS)$ is a projection.

When $n=1$, $W = S_1 T_1^*$.   Note that $T_1^* T_1 = P_{T_1} \in \cP(\cS) = \cQ(\cS)$ and thus $T_1^* T_1 = X_1 X_1^*$ for some $X_1 \in \cS$.  So,  $Q_W = W W^* = S_1 T_1^* T_1 S_1 = S_1 X_1 X_1^* S_1^* = Q_{S_1 X_1} \in \cQ(\cS)$.

Suppose that the result holds for $n \le m$.   Let $V = S_2 T_2^* \cdots S_{m+1} T_{m+1}^*$ and $W = S_1 T_1^* V$.  Then 
\[
WW^* = S_1 T_1^* V V^* T_1 S_1^* = S_1 T_1^* Q_V T_1 S_1^*. \]
Our induction hypothesis implies that $R = Q_V \in \cQ(\cS)$.   From above,  $T_1^* R T_1 \in \cP(\cS) = \cQ(\cS)$, and thus $WW^* = S_1 (T_1^* R T_1) S_1^* \in \cQ(\cS)$ as well.
\end{pf}


\begin{cor} \label{cor2.7}
Suppose that $\cS \subseteq \bofh$ is a semigroup of partial isometries and that $\cP(\cS) \cup \cQ(\cS) \subseteq \cS$.   Then $\cT = \langle \cS \cup \cS^* \rangle$ consists of partial isometries.
\end{cor}

\begin{pf}
The condition that $\cP(\cS) \subseteq \cS$ implies that $\cP(\cS) \subseteq \cQ(\cS)$, since for any $P \in \cP(\cS) \subseteq \cS$, we have that $ P = Q_P \in \cQ(\cS)$.    Similarly, $\cQ(\cS) \subseteq \cP(\cS)$, and so $\cP(\cS) = \cQ(\cS)$.  

The result now follows immediately from Theorem~\ref{thm05}.
\end{pf}


\begin{cor} \label{cor06}
Let $\cS \subseteq \bofh$ be a semigroup of partial isometries and suppose that  the von Neumann algebra $W^*(\cQ(\cS))$ generated by $\cQ(\cS)$ forms a masa in $\bofh$.   Then the semigroup $\cT$ generated by $\cS$ and $\cS^*$ consists of partial isometries.
\end{cor}

\begin{pf}
Simply combine the remarks immediately preceding Theorem~\ref{thm05} with that Theorem.
\end{pf}



\begin{thm} \label{thm07}
Suppose that $\cS \subseteq \bofh$ is a maximal semigroup of partial isometries.   If $\cP(\cS) \cup \cQ(\cS)$ is commutative, then $\cS = \cS^*$.
\end{thm}

\begin{pf}
Let $T \in \cS$.  By Proposition~\ref{prop03}, the semigroup $\cS_0$ generated by $Q_T$ and $\cS$  consists of partial isometries.   Since $\cS$ is a maximal semigroup of partial isometries, it follows that $\cS = \cS_0$.   That is, $Q_T \in \cS$ for all $T \in \cS$ and so $\cQ(\cS) \subseteq \cP(\cS)$.    Applying the same argument to $\cS^*$ (which is also clearly a maximal semigroup of partial isometries) allows us to deduce that $P_T \in \cS$ for all $T \in \cS$, so that $\cP(\cS) \subseteq \cQ(\cS)$.      By Theorem~\ref{thm05}, the semigroup $\cT$ generated by $\cS$ and $\cS^*$ consists of partial isometries.   Once again, the maximality of $\cS$ implies that $\cS = \cT = \cS^*$.
\end{pf}


We remind the reader that Example~\ref{eg04.5} showed that the hypothesis that $\cP(\cS) \cup \cQ(\cS)$ be commutative is not in itself sufficient to guarantee that $\cS$ is contained in a selfadjoint semigroup of partial isometries.   

\smallskip

In proving our next result, we shall need the following result from~\cite{HW1969}.  Recall that an operator $J \in \bofh$ is said to be a \textbf{power partial isometry} if $J^n$ is a partial isometry for all $n \ge 1$.

\begin{thm} \emph{(}\textbf{Halmos-Wallen}\emph{)}  \label{PPI}
Every power partial isometry decomposes as a direct sum whose summands are  unitary operators, pure isometries, pure co-isometries, and truncated shifts.
\end{thm} 


\bigskip

Recall that a if $\cE$ is a collection of projections on a Hilbert space $\hilb$, then a non-zero projection $E \in \cE$ is said to be \textbf{minimal} if $F \in \cE$ and $0 \le F \le E$ implies that $F = 0$ or $F = E$.  The next result generalizes Theorem~3.11 of~\cite{PR2013}.  

\begin{thm} \label{thm08}
Suppose that $\cS$ is an irreducible semigroup of partial isometries and that there exists a minimal non-zero projection $R \in \cP(\cS) \cup \cQ(\cS)$.   Then the semigroup $\cV$ generated by $\cS$ and $\cS^*$ consists of partial isometries.
\end{thm}

\begin{pf}
Considering $\cS^*$, if necessary, we may assume without loss of generality that $R\in\cQ(\cS)$. Suppose that $A\in\cS$ is any partial isometry such that $R=Q_A$. We claim that either $P_A=R$ or $P_AR=0$. Indeed, suppose that $P_AR\ne 0$. Then the operator $A^2$ is not zero. Since the range space of $A^2$ is contained in that of~$A$, the minimality of $R$ yields $Q_{A^2}=Q_A=R$. Decompose $A$ by Theorem~\ref{PPI} into a direct sum of a unitary, pure isometry, pure co-isometry, and truncated shifts. The condition $Q_{A^2}=Q_A$ implies that the pure isometry and the truncated shifts summands are absent. Since the initial space of a pure co-isometry is strictly contained in its final space, again, by the minimality of $R$ we conclude that the co-isometry summand is absent, too. This shows that $A$ is a direct sum of a unitary and zero and proves the claim.

Pick, by Theorem~\ref{irreducible}, an operator $T\in\cS$ such that $ATR\ne 0$. The operator $AT$ is a non-zero member of $\cS$ such that $P_{AT}R\ne 0$. Clearly, $0\ne Q_{AT}\le Q_{A}=R$. Hence, $Q_{AT}=R$ by the minimality of~$R$. Using the claim above with the operator~$AT$, we conclude that $P_{AT}=R$. 

Let $\cR=\{E\in\cP(\cS) \cup \cQ(\cS) \mid E\mbox{ is minimal in }\cP(\cS)\cup \cQ(\cS)$  and there exists  $U\in\cS\mbox{ such that }E=P_U=Q_U\}$. We have shown that $R\in\cR$.  If $E_1,E_2\in\cR$ are two distinct projections, then by minimality we have that $E_1E_2=0$. For each $E\in\cR$, pick a partial isometry $U_E\in\cS$ such that $E=P_{U_E}=Q_{U_E}$.

We claim that $\bigvee_{E\in\cR}E\hilb=\hilb$. Suppose not. Let $P_0=\bigvee_{E\in\cR}E$. Then $P_0^\perp\ne 0$. By Theorem~\ref{irreducible}, there exists $T\in\cS$ such that $P_0^\perp TR\ne 0$. Replacing $T$ with $TU_R\in\cS$, we may assume that the initial space of $T$ is contained in~$R\hilb$. We get:
$$
T = \bordermatrix{ & P_0\hilb & P_0^\perp\hilb \cr
& A & 0 \cr 
& B & 0 }, 
$$
for some operators $A$ and $B$, where $B\ne 0$. If $A\ne 0$, then there exists $E\in\cR$ such that $ET\ne 0$. Consider the operator $U_ET\in\cS$. Clearly, $U_ET\ne 0$ and the initial projection of $U_ET$ is dominated by~$R$. Since $R$ is minimal in $\cP(\cS) \cup \cQ(\cS)$, this projection must coincide with~$R$. It follows that $T$ maps $R\hilb$ isometrically to $E\hilb$, which implies that $B=0$. Since $B$ is assumed to be nonzero, we conclude that $A$ must be equal to zero.

This shows that $T$ is a partial isometry with the initial projection $R$ and the final projection contained in~$P_0^\perp$. Denote $R\hilb$ by~$\hilb_1$, $P_0\hilb\ominus R\hilb$ by~$\hilb_2$ and decompose $P_0^\perp\hilb$ as $\hilb_3\oplus\hilb_4$, where $\hilb_3=T\hilb$ and $\hilb_4=P_0^\perp\hilb\ominus\hilb_3$. Then $T$ can be written as
$$
T = \bordermatrix{ & \hilb_1 & \hilb_2 &  \hilb_3 & \hilb_4 \cr
& 0 & 0 & 0 & 0 \cr 
& 0 & 0 & 0 & 0 \cr 
& C_0 & 0 & 0 & 0 \cr
& 0 & 0 & 0 & 0 \cr},
$$
where $C_0$ is an isometry from $\hilb_1$ onto~$\hilb_3$. Using irreducibility once more, we can find an operator $S\in\cS$ such that $RST\ne 0$. Replacing $S$ with $U_RS$, we get
$$
S = \bordermatrix{ & \hilb_1 & \hilb_2 &  \hilb_3 & \hilb_4 \cr
& A_1 & B_1 & C_1 & D_1 \cr 
& 0 & 0 & 0 & 0 \cr 
& 0 & 0 & 0 & 0 \cr
& 0 & 0 & 0 & 0 \cr},
$$
for some operators $A_1$, $B_1$, $C_1$ and $D_1$, where $C_1\ne 0$. Repeating the arguments used to analyze~$T$, we conclude that $A_1=0$ and $B_1=0$. Consider the operator $ST\in\cS$. The range space of $ST$ is non-zero and contained in $R\hilb$, and therefore must coincide with $R$ by the minimality condition. It follows that the range of $C_1$ is all of~$\hilb_1$. This forces $D_1=0$.

Denote the operator $TS$ by~$V$. We can write
$$
V = \bordermatrix{ & \hilb_1 & \hilb_2 &  \hilb_3 & \hilb_4 \cr
& 0 & 0 & 0 & 0 \cr 
& 0 & 0 & 0 & 0 \cr 
& 0 & 0 & C_0C_1 & 0 \cr
& 0 & 0 & 0 & 0 \cr}.
$$
Clearly, the range space of $V$ is equal to $\hilb_3$ and the initial space of $V$ is contained in~$\hilb_3$. If the initial space of $V$ were not equal to~$\hilb_3$, the fact that $T$ maps $\hilb_1$ isometrically onto $\hilb_3$ would imply that the initial space of $VT$ is non-zero and properly contained in~$\hilb_1$ , contradicting the minimality of~$R$. Therefore, the initial space of $V$ is all of~$\hilb_3$.

Let $F$ be the projection onto $\hilb_3$. Clearly, $F=V^*V=VV^*\in\cP(\cS) \cup \cQ(\cS)$. We claim that $F\in\cR$. To see this, we need to show that $F$ is minimal in $\cP(\cS) \cup \cQ(\cS)$. Indeed, suppose that for some $F_0\in\cP(\cS) \cup \cQ(\cS)$ we have $F_0<F$. If $F_0$ is an initial projection of an operator $T_0\in\cS$, $F_0=T_0^*T_0$, then the initial projection of $T_0T$ is properly contained in~$R$, so that $T_0=0$ and, hence $F_0=0$. If $F_0$ is a final projection of an operator $T_1\in\cS$, $F_0=T_1T_1^*$, then the final projection of $ST_1$ is properly contained in $R$ implying, again, that $F_0=0$. This shows that $F\in\cR$, contradicting the fact that $F\perp P_0$.

We have established that $\bigvee_{E\in\cR}E\hilb=\hilb$. Define $\cS_0$ to be the set of all the partial isometries $W\in\cB(\hilb)$ such that for any pair $E_1,E_2\in\cR$ the operator $E_1WE_2$ is either zero or a partial isometry with $P_{E_1WE_2}=E_2$ and $Q_{E_1WE_2}=E_1$. Let us show that $\cS_0$ is a semigroup.  Let $T,S\in\cS_0$ and $E_1,E_2\in\cR$ be arbitrary. Suppose that $E_1TSE_2\ne 0$. If for every projection $F\in\cR$ the operator $E_1TF$ is zero, then $E_1T=0$ by $\bigvee_{E\in\cR}E\hilb=\hilb$, so that $E_1TSE_2= 0$, a contradiction.   Hence, there exists $F_1\in\cR$ such that $E_1TF_1\ne 0$. Similarly, there exists a projection $F_2\in\cR$ such that $F_2SE_2\ne 0$. It is then easy to see that $E_1T$ is a partial isometry with $P_{E_1T}= F_1$ and $Q_{E_1T}= E_1$. Similarly, $SE_2$ is a partial isometry with $P_{SE_2}= E_2$ and $Q_{SE_2}=F_2$. Since $E_1TSE_2\ne 0$ and any distinct projections in $\cR$ are disjoint, we must have $F_1=F_2$. It follows that $E_1TSE_2$ is a partial isometry with $P_{E_1TSE_2}=E_2$ and $Q_{E_1TSE_2}=E_1$, and $TS\in\cS_0$.

It is obvious that $\cS_0$ is a self-adjoint semigroup. Finally, if $S\in\cS$ and $E_1,E_2\in\cR$, then, by the definition of~$\cR$, $\cS$ has two partial isometries $U_1$ and $U_2$ such that $P_{U_1}=Q_{U_1}=E_1$ and $P_{U_2}=Q_{U_2}=E_2$. The range space of $E_1SE_2$ is equal to the range space of $U_1SU_2$ and is contained in $E_1\hilb$. Similarly, the initial space of $E_1SE_2$ is equal to the initial space of $U_1SU_2$ and is contained in $E_2\hilb$. Since both $E_1$ and $E_2$ are minimal in $\cP(\cS)\cup\cQ(\cS)$ and $U_1SU_2\in\cS$, we conclude that either $E_1SE_2=0$ or $P_{E_1SE_2}=E_2$ and $Q_{E_1SE_2}=E_1$. That is, $\cS\subseteq\cS_0$.
\end{pf}


\begin{rem}\label{brandt}
We notice that the semigroup $\cS_0$ constructed in the proof of Theorem~\ref{thm08} is an example of a well-known Brandt semigroup (see \cite[Chapters 3 and 5]{Howie1995}).
\end{rem}


\begin{rem} \label{rem2.12}
The structure of the semigroup $\cS$ obtained in the proof of Theorem~\ref{thm08} should be compared with the structure of the semigroup $\cS$ which appears in Theorem~3.11 of~\cite{PR2013}.  The finite rank operator $P$ which occurs in the proof of that Theorem is replaced by the minimal projection $R$ appearing in the proof above.

The above proof also shows that every minimal projection in $\cP(\cS) \cup \cQ(\cS)$ actually lies in $\cP(\cS) \cap \cQ(\cS)$.   Furthermore, $\cP(\cS) \cup \cQ(\cS)$ is abelian.   As a consequence of this and the fact that $\vee_{E \in \cR} E \hilb = \hilb$, we find that 
\[
W^*(\cP(\cS) \cup \cQ(\cS)) = W^*(\cP(\cS)) = W^*(\cQ(\cS)). \]
This von Neumann algebra is purely atomic, and as any two minimal projections in $\cP(\cS) \cup \cQ(\cS)$ were shown in the proof to be equivalent, it follows that $W^*(\cP(\cS)\cup \cQ(\cS))$ has uniform multiplicity.   This presages the results of the next section.
\end{rem}


\section{Multiplicity of $W^*(\cQ(\cS))$}

\subsection{}
It will be useful to phrase our next results in the language of direct integrals of Hilbert spaces.   We refer the reader to~\cite{KR1986} for a more detailed treatment than we shall provide.

Recall (see, e.g., \cite[Definition 14.1.6]{KR1986}) that if $\hilb$ is the direct integral of Hilbert spaces over a measure space $(\Omega,\mu)$, 
$$
\hilb=\int_\Omega^\oplus\hilb_p\,d\mu(p),
$$ 
then an operator $T\in\bofh$ is called \textbf{decomposable} if there is a  function $p\mapsto T(p)$ on $\Omega$  such that $T(p)\in\cB(\hilb_p)$ and, for each $x\in\hilb$, $T(p)x(p)=(Tx)(p)$ for almost all~$p$. Decomposable operators will also be denoted as
$$
T=\int_\Omega^\oplus T_p\,d\mu(p).
$$
If, in addition, $T(p)=f(p)I_p$, where $I_p$ is the identity operator on~$\hilb_p$ and $f(p)$ is a scalar, then $T$ is called \textbf{diagonalizable}.


\begin{defn}\label{E_X-notation}
We will use the following notation. Let the Hilbert space be represented as the direct integral $\hilb=\int_\Omega^\oplus\hilb_p\,d\mu(p)$. Let $X$ be a measurable subset of~$\Omega$. We will use the symbol $E_X$ to denote the projection defined by
$$
E_X=\int_X^\oplus I_{p}\,d\mu(p),
$$
where $I_{p}$ is the identity operator on the Hilbert space~$\hilb_p$. These projections will also be referred to as \textbf{standard projections}. The range of a standard projection will be called a \textbf{standard subspace} of~$\hilb$.
\end{defn}


The reader should compare the following result with Corollary~3.12 of~\cite{PR2013}, which states that if $\cS \subseteq \bofh$ is an irreducible semigroup of partial isometries which contains a compact operator, then $\cS$ can be enlarged to a selfadjoint semigroup of partial isometries.

\begin{thm}\label{not-extendable}
There exists an irreducible semigroup $\cS$ of isometries such that $\cP(\cS) \cup \cQ(\cS)$ commutes and $\cS$ cannot be enlarged to a selfadjoint semigroup of partial isometries.
\end{thm}
\begin{pf}
Our example builds on the idea of isometries generating a free semigroup which is weak-operator topology dense in $\bofh$, due to C.~Read~\cite{Rea2005}. Take two isometries with pairwise orthogonal ranges, as exhibited in that paper. The isometries are defined in terms of functions $\phi:[0,\frac{1}{2}]\to[0,1]$ and $\psi:[\frac{1}{2},1]\to[0,1]$ defined by
$$
\phi(t)=2t,\quad\mbox{and}\quad \psi(t)=2t-1,
$$
and specifically chosen functions $\theta_1:[0,\frac{1}{2}]\to\bbC$ and $\theta_2:[\frac{1}{2},1]\to\bbC$ satisfying
$$
\abs{\theta_1 (t)}=\sqrt{2} \mbox{ for all } t\in[0,\tfrac{1}{2}]
$$
and
$$
\abs{\theta_2 (t)}=\sqrt{2} \mbox{ for all } t\in[\tfrac{1}{2},1].
$$
The isometries $T_1,T_2:L_2([0,1])\to L_2([0,1])$ are given by the formulas
$$
(T_1f)(t)=\left\{
\begin{array}{ll}
\theta_1(t)f(\phi(t)),\quad &t\in[0,\frac{1}{2}]\\
0 &t\in[\frac{1}{2},1]
\end{array}\right.
$$
and
$$
(T_2f)(t)=\left\{
\begin{array}{ll}
0 &t\in[0,\frac{1}{2}]\\
\theta_2(t)f(\psi(t)),\quad &t\in[\frac{1}{2},1].
\end{array}\right.
$$
The action of both isometries is just shrinking the support of a given function $f\in L_2([0,1])$ and then multiplying pointwise by a function of modulus~$\sqrt 2$ .

Let $\cK$ be an infinite dimensional, separable Hilbert space and $\cH=L_2([0,1],\cK)$.   Define $S_1,S_2\in\cB(\cH)$ by tensoring $T_1$ and $T_2$ with the identity operator:
$$
(S_1f)(t)=\left\{
\begin{array}{ll}
\theta_1(t)f(\phi(t)),\quad &t\in[0,\frac{1}{2}]\\
0 &t\in[\frac{1}{2},1]
\end{array}\right.
$$
and
$$
(S_2f)(t)=\left\{
\begin{array}{ll}
0 &t\in[0,\frac{1}{2}]\\
\theta_2(t)f(\psi(t)),\quad &t\in[\frac{1}{2},1]
\end{array}\right. .
$$
(The formulas look the same as those for $T_1$ and $T_2$; the values of $f(t)$ in the second set of formulas are vectors in~$\cK$).

Let the symbol $\cU(\cK)$ denote the set of all unitary operators on~$\cK$. Notice that each measurable function $\zeta:[0,1]\to\cU(\cK)$ defines a decomposable unitary operator in~$\cB(\cH)$ by
$$
U_\zeta=\int_{[0,1]}^\oplus\,\zeta(t)\,dt;
$$
this operator does not move the ``fibers'' of $L_2([0,1],\cK)$. Let us denote the family of such unitaries by~$\cU$. 

Consider the semigroup 
$$\cS_0=\langle\{S_1,S_2\}\cup\cU\rangle.$$
Clearly, $\cS_0$ consists of isometries whose ranges are some standard subspaces of~$\cH$.  We claim that $\cS_0$ is irreducible.   Indeed, we may view $\cS_0$ as acting on the tensor product of $L_2[0,1]$ and $\cK$.    Now $\cA := \mathrm{span}\, \langle \{ S_1, S_2\} \rangle$ is weak-operator topology dense in $\bofh$  as shown by Read~\cite{Rea2005}, and it is elementary to see that $\mathrm{span}\ \cU(\cK) = \cB(\cK)$.  Since the algebraic tensor product of two weak-operator topology dense subalgebras of $\cB(L_2([0,1]))$ and $\cB(\cK)$ respectively is weak-operator topology dense in $\cB(L_2 [0,1] \otimes \cK)$, it follows that the algebra generated by $\cS_0$ is weak-operator topology dense in $\bofh$, and thus $\cS_0$ is irreducible.

Let $E=E^*=E^2\in\cB(\cK)$ be a projection such that $\dim E=\dim E^\perp=\infty$. Let $F=I_\cK-E$. Pick two co-isometries onto~$\cK$, $W_1:E\cK\to\cK$ with initial space $E \cK$, and $W_2:F\cK\to\cK$ with initial space $F \cK$. Define a unitary operator $V$ on $\hilb$ by
$$
(Vf)(t)=\left\{
\begin{array}{ll}
\sqrt 2\cdot W_1\Big(E\big(f(\phi(t))\big)\Big), &t\in[0,\frac{1}{2}]\\
\sqrt 2\cdot W_2\Big(F\big(f(\psi(t))\big)\Big),\quad &t\in[\frac{1}{2},1].
\end{array}\right.
$$
That is, $V$ ``splits'' each ``fiber'' into two parts and moves the first part to the interval $[0,\frac{1}{2}]$ and the second part to the interval $[\frac{1}{2},1]$. It is clear that $V$ maps standard subspaces of $\cH$ to standard subspaces.

Now, define the semigroup $\cS\subseteq\cB(\cH)$ by
$$
\cS=\langle\cS_0\cup\{V\}\rangle.
$$
Obviously, $\cS$ consists of isometries. By the properties of $\cS_0$ and $V$, the range space of every isometry from $\cS$ is a standard subspace of~$\cH$. In particular, the set $\cP(\cS) \cup \cQ(\cS)$ is commutative. Finally, $\cS$ is irreducible since it contains an irreducible subsemigroup~$\cS_0$.

Let us show that $\cS$ cannot be enlarged to a self-adjoint semigroup of partial isometries. Indeed, suppose $\cT$ is such an enlargement. Then $\cT$ must contain the projection $Q$ onto the range space of~$S_1$, $L_2([0,\frac{1}{2}],\cK)$. Hence, we must have $QV\in\cT$. The initial projection for $QV$ is the space
$$
P_{QV}=\int_{[0,1]}^{\oplus} E\,dt.
$$
Let $U\in\cB(\cK)$ be a unitary such that $UEU^*$ does not commute with~$E$. By construction, the operator
$$
\widehat U=\int_{[0,1]}^{\oplus} U\,dt
$$
belongs to~$\cS_0$, hence to~$\cT$. It follows that the projection $\widehat UP_{QV}\widehat U^*$ belongs to~$\cT$, too. However, this projection, clearly, does not commute with~$P_{QV}$.
\end{pf}


\begin{rem}\label{infinite-multiplicity}
The semigroup exhibited in the proof of Theorem~\ref{not-extendable} has the property that the set $\cQ(\cS)$ generates an abelian von Neumann algebra of \emph{infinite} uniform multiplicity.  As we shall see, this is crucial.

Our present and main goal is to show that if $\cS \subseteq \bofh$ is an \emph{irreducible} semigroup of partial isometries, then $W^*(\cQ(\cS))$ has uniform multiplicity, and if that multiplicity is finite, then $\cT = \langle \cS \cup \cS^*\rangle$ consists of partial isometries.
\end{rem}


A standard structure theorem for abelian von Neumann algebras $\fN$ (see, for example, Theorem~II.3.5 of~\cite{Dav1996}) shows that if we identify $\fN$ with the algebra $L^\infty(\Omega, \mu)$ for an appropriate compact subset $\Omega \subseteq \bbR$ and positive, regular Borel measure $\mu$ (as in Corollary~II.2.9 of~\cite{Dav1996}), then there exist mutually singular measures $\mu_n << \mu$ such that $\fN$ is unitarily equivalent to the diagonalizable operators on 
\begin{align*}
\hilb 
	&= {\bigoplus}_{1 \le n \le \infty}  \int_{\Omega}^\oplus \hilb_p \, d\mu_n(p) \\
	&= \int_\Omega^\oplus \hilb_p \, d\mu(p),
\end{align*}
where $\dim\, \hilb_p = n$ for almost every $p$ in the support of $\mu_n$.

It is with respect to this decomposition that we will be decomposing the abelian von Neumann algebra $W^*(\cQ(\cS))$ below.

\bigskip
	
We remind the reader that 
\[
E_X=\int_X^\oplus I_{p}\,d\mu(p), \]
where $I_{p}$ is the identity operator on the Hilbert space~$\hilb_p$.


\begin{prop}\label{set-mapping}
Let $\cS \subseteq \bofh$ be a semigroup of partial isometries such that $\cQ(\cS)$ is commutative. Let $W^*(\cQ(\cS))$ be represented as the  multiplication operators on $\hilb=\int_\Omega^\oplus\hilb_p\,d\mu(p)$. If $T\in\cS$, then for every measurable set $X\subseteq\Omega$, the operator $TE_X$ is a partial isometry and there is a measurable set $Y\subseteq\Omega$ such that 
$$
Q_{TE_X}=E_Y.
$$
Moreover, if $(X_n)_n$ is a sequence of measurable subsets of $\Omega$ and $Y_n\subseteq\Omega$ are such that $Q_{TE_{X_n}}=E_{Y_n}$ for all $n \ge 1$, then 
$$
Q_{TE_{\cup_{n=1}^\infty X_n}}=E_{\cup_{n=1}^\infty Y_n}.
$$
\end{prop}

\begin{pf}
The fact that $TE_X$ is a partial isometry for every measurable set $X\subseteq\Omega$ follows immediately from Theorem~\ref{HalmosWallen}. We need to prove the assertion about the~$Q_{TE_X}$.

Denote the von Neumann algebra $W^*(\cQ(\cS))$ by~$\fM$. Notice that if $T\in\cS$, then the final projection $Q_T$ is of the form $E_{Y_T}$, for some measurable subset $Y_T$ of~$\Omega$. Moreover, the $\sigma$-algebra of measurable subsets of $\Omega$ is equal to the completion of the $\sigma$-algebra generated by the set $\left\{Y_T\mid T\in\cS\right\}$.


Define a set $\Lambda$ of measurable subsets of $\Omega$ by
\begin{align*}
\Lambda=\{X\subseteq\Omega\mid\ X\mbox{ is measurable }&\mbox{and }Q_{TE_X}\in\fM\}.
\end{align*}
Our goal is to prove that $\Lambda$ coincides with the set of all measurable subsets of~$\Omega$.

Notice that if $X$ is a measurable subset such that $Q_{TE_X}\in\fM$, then $Q_{TE_X}$ must be of the form $E_Y$ for some measurable subset $Y$ of~$\Omega$. Also, observe that by Proposition~\ref{prop03}, $\Lambda$ contains every measurable subset $X\subseteq\Omega$ such that $E_X=Q_S$ for some operator $S\in\cS$. It is also clear that $\Omega\in\Lambda$ and $\varnothing \in\Lambda$.

Let us first establish that $\Lambda$ forms a Boolean algebra. Let $Y_0\subseteq\Omega$ be such that $Q_T=E_{Y_0}$. Fix an arbitrary $X\in\Lambda$. By the definition of~$\Lambda$, there is a measurable subset $Y$ of $\Omega$ such that $Q_{TE_X}=E_{Y}$. The projection $Q_{TE_X}$ is a subprojection of $Q_T$,  hence $Y\subseteq Y_0$. Therefore
$$
Q_{TE_{\Omega\setminus X}}=Q_{T(I-E_X)}=T(I-E_X)T^*=TT^*-TE_XT^*=Q_T-Q_{TE_X}=E_{Y_0\setminus Y}.
$$
So, $\Lambda$ is stable under taking complements. Let $X_1, X_2\in\Lambda$ be arbitrary. There exist two measurable sets $Y_1$ and $Y_2$ such that $Q_{TE_{X_1}}=E_{Y_1}$ and $Q_{TE_{X_2}}=E_{Y_2}$. By Lemma~\ref{lem01alternate}(a) used with the (commutative) semigroup of all the projections in~$\fM$, we get
$$
Q_{TE_{X_1\cap X_2}}=TE_{X_1}E_{X_2}T^*=Q_{TE_{X_1}}Q_{TE_{X_2}}=E_{Y_1}E_{Y_2}=E_{Y_1\cap Y_2}.
$$
This shows that $\Lambda$ forms a Boolean algebra. 

Let us now show that $\Lambda$ is a $\sigma$-algebra. Suppose that $\{X_n\}_{n\in\bbN}$ is an increasing sequence of sets in~$\Lambda$ and let $Z=\cup_{n\in\bbN}X_n$. Clearly, $Z\subseteq\Omega$ is measurable. We claim that $Z\in\Lambda$. Indeed, it follows from Theorem~\ref{HalmosWallen} that $TE_Z$ is a partial isometry. For each $n\in\bbN$, let $Y_n\subseteq\Omega$ be such that $Q_{TE_{X_n}}=E_{Y_n}$. Let $Y=\cup_{n\in\bbN}Y_n$. It is clear that $TE_Z$ is an extension of every $TE_{X_n}$ ($n\in\bbN$). Hence, $Q_{TE_Z}$ dominates each~$E_{Y_n}$. Thus
$$
Q_{TE_Z}\ge E_Y.
$$ 
On the other hand, suppose that $h\in E_Y^\perp\hilb$. Clearly, $h\perp E_{Y_n}\hilb$ for every $n\in\bbN$. That is, for every $f\in\hilb$ and $n\in\bbN$ we have $h\perp TE_{X_n}f$. Notice that $TE_Z=\lim_{n\to\infty}TE_{X_n}$ in the strong operator topology. We conclude that $h\perp TE_Zf$ for every $f\in\hilb$. This shows that
$$
Q_{TE_Z}\le E_Y,
$$
so that 
$$
Q_{TE_Z}=E_Y.
$$
We established that if $\{X_n\}_{n\in\bbN}$ is an increasing sequence of sets in~$\Lambda$, then $\cup_{n\in\bbN} X_n\in\Lambda$ and
$$
Q_{TE_{\cup X_n}}=E_{\cup Y_n},
$$
where $X_n$ is such that $Q_{TE_{X_n}}=E_{Y_n}$. Similarly, one can show that if $\{X_n'\}_{n\in\bbN}$ is a decreasing sequence of sets in~$\Lambda$, then $\cap_{n\in\bbN} X'_n\in\Lambda$ and
$$
Q_{TE_{\cap Y'_n}}=E_{\cap X'_n},
$$
where $Y'_n$ is such that $Q_{TE_{X'_n}}=E_{Y'_n}$.

By the Monotone Class Theorem (see, e.g., \cite[Theorem 1.9.3(i)]{Bog2007}), $\Lambda$ contains the $\sigma$-algebra generated by the set $\left\{X\subseteq\Omega\mid E_X=Q_S\mbox{ for some }S\in\cS\right\}$. Since $\Lambda$ also has the property that if $X\in\Lambda$ and $X_0\subseteq\Omega$ is such that $\mu(X_0)=0$, then $X\cup X_0\in\Lambda$, it follows that $\Lambda$ contains every measurable subset of~$\Omega$.

The ``moreover'' part is clear from the construction.
\end{pf}

The following Corollary extends Proposition~\ref{prop03}.
\smallskip

\begin{cor}\label{W*-enrich}
Let $\cS$ be a semigroup of partial isometries such that $\cQ(\cS)$ is commutative. If $E\in W^*(\cQ(\cS))$ is an arbitrary projection, then the semigroup $\cT=\langle\cS \cup \{E\}\rangle$ consists of partial isometries and has the property that $W^*(\cQ(\cT))=W^*(\cQ(\cS))$.
\end{cor}

\begin{pf}
Since we have assumed that $\cS$ is untial, an arbitrary member of $\cT$ has the form
$$
S=T_1ET_2E\dots T_{n-1}ET_n,
$$
with $T_i\in\cS$. We will show by induction that $S$ is a partial isometry and $Q_S\in W^*(\cQ(\cS))$. If $n=1$, the statement is trivial. Suppose that the statement is valid for~$n$; we need to establish it for $n+1$. Write $S$ as
$$
S=T_1ES_0,
$$
where $S_0$ is, by the induction hypothesis, a partial isometry such that $Q_{S_0}\in W^*(\cQ(\cS))$. Notice that every projection in $W^*(\cQ(\cS))$ is of the form $E_X$, for some measurable subset $X$ of $\Omega$.  Fix measurable sets $X$ and $X_0$ such that $E=E_X$ and $Q_{S_0}=E_{X_0}$. By Theorem~\ref{HalmosWallen}, the operator $ES_0$ is a partial isometry. Moreover, $Q_{ES_0}=EQ_{S_0}E^* = E Q_{S_0} E = E_{X\cap X_0}$. Since the initial projection of $T_1$ commutes with $W^*(\cQ(\cS))$, it follows, again, by Theorem~\ref{HalmosWallen}, that $T_1ES_0$ is a partial isometry. Finally, by Lemma~\ref{lem01alternate}, $Q_{T_1ES_0}=T_1EQ_{S_0}ET_1^*=T_1E_{X\cap X_0}T_1^*=Q_{T_1E_{X\cap X_0}}$ which, by Proposition~\ref{set-mapping}, belongs to~$W^*(\cQ(\cS))$.
\end{pf}


\begin{rem}\label{S-max-rich}
It follows from Corollary~\ref{W*-enrich} that the semigroup $\cS_{max}$ constructed in Proposition~\ref{prop04} contains all the projections from $W^*(\cQ(\cS))=W^*(\cQ(\cS_{max}))$.
\end{rem}


\begin{lem}\label{decomposable}
Let the Hilbert space $\hilb$ be represented as a direct integral, $\hilb=\int_{\Omega}^\oplus \hilb_p\,d\mu(p)$. If $T\in\bofh$ is such that $E_X$ is $T$-invariant for every measurable subset $X$ of~$\Omega$, then $T$ is decomposable.
\end{lem}

\begin{pf}
Let $X\subseteq \Omega$ be any measurable subset of $X$. Then $I - E_X = E_{\Omega \setminus X}$.   Our hypotheses show that both $E_X \hilb$ and $(E_X \hilb)^\perp = E_{\Omega\setminus X}\hilb$ are invariant for $T$, and thus $E_X \hilb$ is orthogonally reducing for $T$ for all $X \subseteq \Omega$ measurable.  

It follows that $T E_X = E_X T$ for all measurable $X \subseteq \Omega$.  The lemma now follows from \cite[Theorem 14.1.10]{KR1986} and the Double Commutant Theorem.
\end{pf}


It was established in~\cite{PR2013} that operators in a self-adjoint semigroup of partial isometries can be simultaneously represented as generalized composition operators. The proof of the following lemma shows that a similar representation is possible for operators in non-self-adjoint semigroups of partial isometries, at least on some parts of their domains. 

Notice that if $\cS$ is a semigroup of partial isometries with commuting set of final projections $\cQ(\cS)$ and $T\in\cS$, then the initial projection $P_T$ of $T$ belongs to $W^*(\cQ(\cS))'$ by Theorem~\ref{HalmosWallen}. Hence, if we represent $W^*(\cQ(\cS))$ as the set of multiplication operators on $\hilb=\int_\Omega^\oplus\hilb_p\,d\mu(p)$, the operator $P_T$ becomes a decomposable operator with respect to this decomposition. That is, we can write 
$$
P_T=\int_\Omega^\oplus P_p\,d\mu(p).
$$
Let us now state the lemma which will play a key role for the remainder of the paper.

\begin{lem}\label{composition-representation}
Let $\cS$ be a semigroup of partial isometries acting on a separable
Hilbert space such that $\cQ(\cS)$ is commutative and let $T\in\cS$ be arbitrary. Represent $W^*(\cQ(\cS))$ as the diagonalizable operators on $\hilb=\int_\Omega^\oplus\hilb_p\,d\mu(p)$ and write $P_T=\int_\Omega^\oplus P_p\,d\mu(p)$. Let $1 \le k < \infty$ be an integer and let $X\subseteq\Omega$ be a measurable set such that $\rk P_p=k$ for every $p\in X$. Let $Y\subseteq\Omega$ be  a measurable set satisfying  $Q_{TE_X}=E_Y$ (such exists by Proposition~\ref{set-mapping}). If $\dim\hilb_p\ge k$ for all $p\in Y$, then $\dim \hilb_p=k$ for almost all $p\in Y$.
\end{lem}

\begin{pf}
We will prove more: we will show that, perhaps after changing $X$ by a set of measure zero, there exist a measurable injective map $\phi:X\to\Omega$ and a family of operators $A_T(p):P_p(\hilb_p)\to\hilb_{\phi(p)}$ ($p\in X$) such that $Tf\in\int_{\phi(X)}^\oplus\hilb_p\,d\mu(p)$ for every $f\in\int_{X}^\oplus P_p(\hilb_p)\,d\mu(p)$ and, moreover,
\begin{equation}
\label{eq:2}
(Tf)(\phi(p))=A_T(p)f(p),\quad p\in X.
\end{equation}
Since the final space of the operator $TE_X$ is of the form $E_Y$, the range of each $A_p$ must be all of $\hilb_{\phi(p)}$, so that formula~\eqref{eq:2} implies the conclusion of the lemma.

The proof is by induction on~$k$. Suppose first that $k=1$.

Denote the $\sigma$-algebra of the measurable subsets of $X$ by~$\Sigma$. By Proposition~\ref{set-mapping}, for every $Z\in\Sigma$, there is a measurable subset $V$ of $\Omega$ such that $Q_{TE_Z}=E_V$. Put
\begin{equation}
\label{eq:3}
\Lambda=\{V\subseteq Y\mid V\mbox{ is measurable and }Q_{TE_Z}=E_V\mbox{ for some }Z\in\Sigma\}.
\end{equation}
We will show that $\Lambda$ coincides with the set of all measurable subsets of~$Y$. 

Let $V\subseteq Y$ be an arbitrary measurable subset. Consider the operator $S=E_VTE_{X}$. By Corollary~\ref{W*-enrich}, this is a partial isometry, and its initial projection is
$$
P_S=(E_VTE_{X})^*E_VTE_{X}=E_{X}T^*E_VTE_{X}\le E_{X}T^*TE_{X}=E_{X}P_TE_{X}=P_TE_{X}.
$$
However,
$$
P_TE_{X}=\int_{X}^\oplus P_p\,d\mu(p),
$$
and for each $p\in X$ the rank of $P_p$ is equal to one. It follows that there exists a measurable subset $Z$ of $X$ such that
\begin{equation}
\label{eq:4}
P_S=\int_Z^\oplus P_p\,d\mu(p)=P_TE_Z.
\end{equation}
Then
\begin{align*}
Q_{TE_Z} &=TE_ZT^*=TP_TE_ZT^*=TP_ST^*=TS^*ST^*=\\
&=T(E_VTE_{X})^*E_VTE_{X}T^*=TE_{X}T^*E_VTE_{X}T^*=\\
&=Q_{TE_{X}}E_VQ_{TE_{X}}=Q_{TE_{X}}E_V=E_{Y}E_V=E_V.
\end{align*}	 
This shows that $V\in\Lambda$.

Define a measure $\nu$ on the measurable space $(Y,\Lambda)$ by letting
$$
\nu(V)=\mu(Z),\quad\mbox{where }Z \subseteq X \mbox{ is such that }Q_{TE_Z}=E_V.
$$
Firstly, let us see that $\nu$ is well-defined. Suppose that $Z_1, Z_2 \subseteq X$ are measurable sets and that $Q_{TE_{Z_1}}=Q_{TE_{Z_2}}$. Denote by $Z$ the symmetric difference of $Z_1$ and $Z_2$, $Z=Z_1\Delta Z_2$. It is easy to see that $Q_{TE_Z}=0$. This means that $TE_Z=0$. However, the fact that $\rk(P_p)=k>0$ for all $p\in X$ implies that $\mu(Z)=0$, so that $\mu(Z_1)=\mu(Z_2)$.

Secondly, the fact that $\nu$ is a measure follows from the ``moreover'' part of Proposition~\ref{set-mapping}.

It is not hard to see that the measures $\mu$ and $\nu$ are absolutely continuous with respect to each other on the set~$Y$. In particular, $\nu$ is a regular Borel measure.

Recall (see~\cite[\S9.3]{Bog2007}) that a measure algebra of a measure space $(\Omega,\Sigma,\mu)$ is the metric space of all the equivalence classes of elements of~$\Sigma$ endowed with the metric $d(A,B)=\mu(A\Delta B)$. It follows that the measure algebras corresponding to the measure spaces $(X,\Sigma,\mu)$ and $(Y,\Lambda,\nu)$ are isomorphic (in the sense of \cite[Definition 9.3.1]{Bog2007}). By \cite[Corollary 9.5.2]{Bog2007}, there exist two sets $N\subseteq X$ and $M\subseteq Y$ such that $\mu(N)=\nu(M)=0$ and a bijective map $\phi:X\setminus N\to Y\setminus M$ such that $\phi(\Sigma_{X\setminus N})=\Lambda_{Y\setminus M}$ and $\nu(\phi(Z))=\mu(Z)$ for all $Z\in\Sigma_{X\setminus N}$. Relabeling, we may assume that $N=M=\emptyset$.

Pick an operator $U:E_Y\hilb\to P_TE_X\hilb$ defined on each $f\in E_Y\hilb$ by the formula
$$
(Uf)(p)=U_p\big(f(\phi(p))\big)\cdot \sqrt{ \frac{d\nu}{d\mu}(\phi(p))},
$$
where each $U_p$ is a surjective partial isometry $\hilb_{\phi(p)}\to P_p\hilb_p$, the function $p\mapsto U_p$ is measurable, and $\frac{d\nu}{d\mu}$ denotes the Radon-Nikodym derivative. It is not hard to see that $U$ is a partial isometry. Define the operator
$$
\widetilde T = TU :E_Y \hilb \to E_Y \hilb.
$$
Since $Q_{TE_Z}=E_{\phi(Z)}$ for every $Z\in\Sigma$, it is easy to see that $\widetilde T$ leaves every space of the form $E_{\phi(Z)}$ invariant. It follows from Lemma~\ref{decomposable} that $\widetilde T$ is a decomposable operator on $\int_Y^\oplus \hilb_p d\mu(p)$. That is, 
$$
(\widetilde Tf)(\phi(p))=T_p\big(f(\phi(p))\big),\quad f\in E_Y\hilb,\ p\in X,
$$
where $T_p\in\cB(\hilb_{\phi(p)})$. Since $T|_{E_X \hilb} =\widetilde TU^*$, we get, for every $f\in P_TE_X\hilb$, 
$$
(Tf)(\phi(p))=T_p\big((U^*f)(\phi(p))\big)=T_pU^*_p(f(p))\cdot \sqrt{\frac{d\nu}{d\mu}(\phi(p))},\quad p\in X.
$$
This verifies the formula~\eqref{eq:2} for $k=1$.

Suppose now that the statement has been proved for all $n\le k$. Let us establish it for~$k+1$. The proof is analogous to the case $k=1$. Again, fix a measurable set $Y$ such that $Q_{ET_X}=E_Y$ and define the set $\Lambda$ as in formula~\eqref{eq:3}. Let us show that $\Lambda$ consists of all the measurable subsets of~$Y$. As in the case $k=1$, pick an arbitrary measurable subset $V$ of $Y$ and define $S=E_VTE_X$. Again, it is easy to see that 
$$
P_S\le P_TE_{X}=\int_{X}^\oplus P_p\,d\mu(p).
$$
It follows that
$$
P_S=\int_{X}^\oplus P'_p\,d\mu(p),
$$
where $P'_p$ is a projection on $\hilb_p$ such that $P'_p\le P_p$. For each $n\le k+1$, let $X_n=\{p\in X\mid\rk(P'_p)=n\}$. For each $0<n\le k$, consider the operator 
$$
S_n=SE_{X_n}.
$$
By Corollary~\ref{W*-enrich}, the semigroup $\cV_n=\langle\cS,E_X,E_V,E_{X_n}\rangle$ consists of partial isometries and its final projections generate the same von Neumann algebra as those of~$\cS$. Also, $S\in\cV_n$ and $S_n\in\cV_n$. By Proposition~\ref{set-mapping}, there is a measurable set $Y_n\subseteq Y$ such that $Q_{S_n}=Q_{SE_{X_n}}=E_{Y_n}$. Hence, by the induction hypothesis, for almost every $p\in Y_n$, $n \le k$, we must have $n\ge\dim\hilb_p$. However, $Y_n\subseteq Y$, and therefore $\dim\hilb_p\ge (k+1)>n$ for all $p\in Y$. It follows that $\mu(Y_n)=0$, which implies $\mu(X_n)=0$. 

This shows that there exists a measurable set $Z\subseteq X$ such that $P_S=P_TE_Z$ (take $Z=X_{k+1}$). This is the same conclusion as obtained in formula~\eqref{eq:4} corresponding to the case $k=1$. The rest of the proof repeats the proof of the case $k=1$ verbatim.
\end{pf}


The next two theorems are simple corollaries of Lemma~\ref{composition-representation}.

\begin{thm}\label{uniform-multiplicity}
Let $\cS$ be a semigroup of partial isometries such that $\cQ(\cS)$ is commutative. If $\cS$ is irreducible, then $W^*(\cQ(\cS))$ has uniform multiplicity.
\end{thm}

\begin{pf}
As before, represent $W^*(\cQ(\cS))$ as the multiplication operators on
$$
\hilb=\int_{\Omega}^\oplus\hilb_p\,d\mu(p).
$$
Write 
$$
\Omega=\left(\cup_{n=1}^\infty\Omega_n\right)\cup\Omega_{\infty},
$$
where $\Omega_n=\{p\in\Omega\mid\dim\hilb_p=n\}$, $n\in\bbN\cup\{\infty\}$. If $W^*(\cQ(\cS))$ does not have a uniform multiplicity, then there exist $n\in\bbN$ and $m\in\bbN\cup\{\infty\}$ such that $n<m$ and $\mu(\Omega_n)>0$ and $\mu(\Omega_m)>0$. Using Corollary~\ref{W*-enrich}, we may assume that $E_{\Omega_n}\in\cS$ and $E_{\Omega_m}\in\cS$.

By irreducibility and Theorem~\ref{irreducible}, there exists an operator $T\in\cS$ such that $E_{\Omega_m}TE_{\Omega_n}\ne 0$. Since $E_{\Omega_m}TE_{\Omega_n}\in\cS$, this clearly contradicts the conclusion of Lemma~\ref{composition-representation}.
\end{pf}


\begin{thm}\label{finite-multiplicity}
Let $\cS$ be a semigroup of partial isometries such that $W^*(\cQ(\cS))$ has uniform finite multiplicity. If $T\in\cS$ then $P_T\in W^*(\cQ(\cS))$ for all $T\in\cS$.
\end{thm}

\begin{pf}
Again, represent $W^*(\cQ(\cS))$ as the multiplication operators on
$$
\hilb=\int_{\Omega}^\oplus\hilb_p\,d\mu(p).
$$
Denote the multiplicity of $W^*(\cQ(\cS))$ by~N. That is, $\dim\hilb_p=N$ for all $p\in\Omega$. Since $P_T\in(W^*(\cQ(\cS)))'$ by Theorem~\ref{HalmosWallen}, $P_T$ is a decomposable operator,
$$
P_T=\int_\Omega^\oplus P_p\,d\mu(p).
$$
Suppose that $P_T\not\in W^*(\cQ(\cS))$. Then $P_T$ is not diagonalizable. For $n\le N$, let $\Omega_n=\{p\in\Omega\mid\dim(\hilb_p)=n\}$. Then each $\Omega_n$ is a measurable set. It follows that at least one of the sets $\Omega_n$, $1\le n<N$, has non-zero measure. This, however, contradicts the conclusion of Lemma~\ref{composition-representation}.
\end{pf}


\begin{cor}\label{masa-ext}
Let $\cS$ be a semigroup of partial isometries such that $\cQ(\cS)$ is commutative and $W^*(\cQ(\cS))$ has uniform finite multiplicity. If $T\in\cS$ then the semigroup $\cS_0$ generated by $Q_T$ and by $\cS$ consists of partial isometries and satisfies the condition $\cP(\cS_0)\cup\cQ(\cS_0)\subseteq W^*(\cQ(\cS))$.
\end{cor}

\begin{pf}
The fact that $\cS_0$ consists of partial isometries was proved in Proposition~\ref{prop03}, so we only need to prove the statement about the projections. Clearly, $W^*(\cQ(\cS))\subseteq W^*(\cQ(\cS_0))$. By Proposition~\ref{prop03}, $\cQ(\cS_0)\subseteq W^*(\cQ(\cS))$, so that $W^*(\cQ(\cS_0))=W^*(\cQ(\cS))$. In particular, $W^*(\cQ(\cS_0))$ has uniform finite multiplicity. It readily follows from Theorem~\ref{finite-multiplicity} that $\cP(\cS_0)\subseteq W^*(\cQ(\cS_0))$, so that $\cP(\cS_0)\cup\cQ(\cS_0)\subseteq W^*(\cQ(\cS_0))=W^*(\cQ(\cS))$.
\end{pf}


\begin{prop} \label{max-prop}
Let $\cS$ be a semigroup of partial isometries such that $\cQ(\cS)$ is commutative and $W^*(\cQ(\cS))$ has uniform finite multiplicity.   
Then 
\begin{enumerate}
	\item[(a)]
	there exists a semigroup $\cS_{\textrm{max}}$ of partial isometries which is maximal with respect to the following conditions:
	\begin{enumerate}
		\item[(i)]
		$\cS \subseteq \cS_{\textrm{max}}$, and 
		\item[(ii)]
		$\cP(\cS_{max}) \cup \cQ(\cS_{max}) \subseteq W^*(\cQ(\cS))$.
	\end{enumerate}
	\item[(b)]
	Furthermore, $\cP(\cS_{max}) \cup \cQ(\cS_{max}) \subseteq \cS_{max}$.
\end{enumerate}
\end{prop}

\begin{pf}

\begin{enumerate}

	\item[(a)]
	Let 
	\begin{align*}
	\fS &:=  \{ \cR \subseteq \bofh: \cR \mbox{ is a semigroup of partial isometries  }   \\
		& \ \ \ \ \ \ \ \ \ \ \ \ \ \ \ \ \ \ \ \ 
			\ \ \ \ \ \ \ \mbox{ containing } \cS, \mbox{ and } \cP(\cR) \cup \cQ(\cR) \subseteq W^*(\cQ(\cS)) \},
	\end{align*}	
	partially ordered  by inclusion.  A standard application of Zorn's Lemma shows that $\fS$ admits a maximal element, $\cS_{\textrm{max}}$, which clearly satisfies conditions (i) and (ii).
	\item[(b)]

	Suppose that $T \in \cS_{\textrm{max}}$.   Then Corollary~\ref{masa-ext} combined with the maximality of $\cS_{\textrm{max}}$ implies 
	that $Q_T \in \cS_{\textrm{max}}$.    Next, we show that $\cP(\cS_{\textrm{max}}) \subseteq \cS_{\textrm{max}}$.  Before doing so, it is worth noting that  $W^*(\cP(\cS_{\textrm{max}}))= W^*(\cQ(\cS))$.
	
	Indeed, $\cP(\cS_{\textrm{max}}) \supseteq \cQ(\cS_{\textrm{max}})$, since $Q \in \cQ(\cS_{\textrm{max}})$ implies $Q \in \cS_{\textrm{max}}$ whence $P_Q = Q \in \cP(\cS_{\textrm{max}})$.  Hence $\cP(\cS_{\textrm{max}}) \supseteq \cQ(\cS)$ and so $W^*(\cP(\cS_{\textrm{max}})) \supseteq W^*(\cQ(\cS))$.
	
	Conversely, $\cP(\cS_{\textrm{max}}) \subseteq W^*(\cQ(\cS))$ implies $W^*(\cP(\cS_{\textrm{max}})) \subseteq W^*(\cQ(\cS))$, and so these last two algebras are in fact equal.

	Let $\cV := \cS_{\textrm{max}}^*$.   Then $\cV$ is a semigroup of partial isometries, $\cQ(\cV) = \cP(\cS_{\textrm{max}})$ is commutative and $W^*(\cQ(\cV)) = W^*(\cP(\cS_{\textrm{max}}))$ has uniform multiplicity.  Applying part (a) to $\cV$ yields a semigroup $\cV_{\textrm{max}}$ which contains $\cV$ and for which 
	\[ 
	\cP(\cV_{\textrm{max}}) \cup \cQ(\cV_{\textrm{max}}) \subseteq W^*(\cQ(\cV)) = W^*(\cP(\cS_{\textrm{max}})) \subseteq W^*(\cQ(\cS)). \]
	
Now, $\cS_0 := \cV_\textrm{max}^* \supseteq \cV^* = \cS_{\textrm{max}}$ and 
	\[
	\cP(\cS_0) \cup \cQ(\cS_0) = \cQ(\cV_{\textrm{max}}) \cup \cP(\cV_{\textrm{max}}) \subseteq W^*(\cQ(\cS)). \]
	The maximality of $\cS_{\textrm{max}}$ therefore implies that $\cS_0 = \cS_{\textrm{max}}$, and so $\cV_{\textrm{max}} = \cV$.  
	
	From the first paragraph of part (b), we see that $\cQ(\cV_{\textrm{max}}) \subseteq \cV_{\textrm{max}}$, i.e. $\cQ(\cV) \subseteq \cV$.    Hence $\cP(\cS_{\textrm{max}}) \subseteq \cV$, whence $\cP(\cS_{\textrm{max}}) = \cP(\cS_{\textrm{max}})^* \subseteq \cV^* = \cS_{\textrm{max}}$.
	
	\end{enumerate}
	\end{pf}

The following is one of the main results of the paper, and generalizes Corollary~\ref{cor06}.

\begin{thm} \label{thm2.24}
Let $\cS$ be a semigroup of partial isometries such that $\cQ(\cS)$ is commutative and $W^*(\cQ(\cS))$ has uniform finite multiplicity.     Then the semigroup $\cT = \langle \cS \cup \cS^* \rangle$ consists of partial isometries.  
\end{thm}

\begin{pf}
It follows from Proposition~\ref{max-prop} above that we can embed $\cS$ in a semigroup $\cS_{\textrm{max}}$ of partial isometries for which $\cP(\cS_{\textrm{max}}) \cup \cQ(\cS_{\textrm{max}}) \subseteq \cS_{\textrm{max}}$.  By Corollary~\ref{cor2.7}, the semigroup $\cT_{\textrm{max}} := \langle \cS_{\textrm{max}} \cup \cS^*_{\textrm{max}} \rangle$ consists of partial isometries.  Since $\cT \subseteq \cT_{\textrm{max}}$, we are done.
\end{pf}


\begin{rem} \label{rem3.15}
If $\cS \subseteq \bofh$ is an irreducible semigroup of partial isometries and if $\cQ(\cS)$ is commutative, then $W^*(\cQ(\cS))$ has uniform multiplicity by Theorem~\ref{uniform-multiplicity}.  

If that multiplicity is finite, then $\cT = \langle \cS \cup \cS^* \rangle$ consists of partial isometries by Theorem~\ref{thm2.24}.   

If that multiplicity is infinite, then Theorem~\ref{not-extendable} shows that $\cT$ need not consist of partial isometries.  On the other hand,   let $\cU(\hilb)$ denote the set of unitary operators acting on $\hilb$ (infinite-dimensional, complex and separable) with orthonormal basis $\{ e_n \}_{n=1}^\infty$, and given $x, y \in \hilb$, denote by $x \otimes y^*$ the rank-one operator $x\otimes y^*(z) = \langle z, y \rangle\, x$ for all $z \in \hilb$.  If we let $\cS \subseteq \cB(\hilb \otimes \hilb)$ be the irreducible semigroup 
\[
\cS := \langle \{U \otimes (e_i \otimes e_j^*): U \in \cU(\hilb), 1 \le i, j \le \infty \} \rangle, \]
then it is reasonably straightforward to check that $W^*(\cQ(\cS))$ has uniform infinite multiplicity and that $\cS = \cS^*$.
\end{rem}















\bibliographystyle{plain}
\bibliography{BMPR2013}

\begin{thebibliography}{10}

\bibitem{Bar1976}
B.A. Barnes.
\newblock Representations of the $\ell_1$-algebra of an inverse semigroup.
\newblock {\em Trans. Amer. Math. Soc.}, 218:361--396, 1976.

\bibitem{Bog2007}
V.I. Bogachev.
\newblock {\em Measure theory. Vol. I, II.}
\newblock Springer-Verlag, Berlin, 2007.

\bibitem{Cun1977}
J.~Cuntz.
\newblock Simple ${C}^*$-algebras generated by isometries.
\newblock {\em Commun. Math. Phys.}, 57:173--185, 1977.

\bibitem{CK1980}
J.~Cuntz and W.~Krieger.
\newblock A class of ${C}^*$-algebras and topological {M}arkov chains.
\newblock {\em Inventiones math.}, 56:251--268, 1980.

\bibitem{Dav1996}
K.R. Davidson.
\newblock {\em ${C}^*$-algebras by example}, volume~6 of {\em Fields Institute
  Monographs}.
\newblock Amer. Math. Soc., Providence, RI, 1996.

\bibitem{DKP2001}
K.R. Davidson, E.~Katsoulis, and D.R. Pitts.
\newblock The structure of free semigroup algebras.
\newblock {\em J. Reine Angew. Math.}, 533:99--125, 2001.

\bibitem{DP1985}
J.~Duncan and A.L.T. Paterson.
\newblock ${C}^*$-algebras of inverse semigroups.
\newblock {\em Proc. Edinburgh Math. Soc.}, 28:41--58, 1985.

\bibitem{HW1969}
P.R. Halmos and L.J. Wallen.
\newblock Powers of partial isometries.
\newblock {\em J. Math. Mech.}, 19:657--663, 1969/1970.

\bibitem{Howie1995}
J.M. Howie.
\newblock {\em Fundamentals of Semigroup Theory}.
\newblock London Mathematical Society Monographs New Series. Clarendon Press,
  Oxford, 1995.

\bibitem{KR1986}
R.~Kadison and J.~Ringrose.
\newblock {\em Fundamentals of the theory of operator algebras. Vol. II.
  Advanced theory.}, volume 100 of {\em Pure and Applied Mathematics}.
\newblock Academic Press, Inc., Orlando, FL, 1986.

\bibitem{Pet1984}
M.~Petrich.
\newblock {\em Inverse Semigroups}.
\newblock Pure and Applied Mathematics. John Wiley and Sons, New York, 1984.

\bibitem{Pop1996}
G.~Popescu.
\newblock Non-commutative disc algebras and their representations.
\newblock {\em Proc. Amer. Math. Soc.}, 124:2137--2148, 1996.

\bibitem{PR2013}
A.I. Popov and H.~Radjavi.
\newblock Semigroups of partial isometries.
\newblock {\em to appear in Semigroup Forum}, 2013; Corrigendum: Semigroups of
  partial isometries, \emph{to appear in Semigroup Forum}, 2013.

\bibitem{Pre1954}
G.B. Preston.
\newblock Representation of inverse semigroups.
\newblock {\em J. London Math. Soc.}, 29:411--419, 1954.

\bibitem{Rae2005}
I.~Raeburn.
\newblock {\em Graph algebras}, volume 103 of {\em CBMS Regional Conference
  Series in Mathematics}.
\newblock Amer. Math. Soc., Washington, D.C., 2005.

\bibitem{Rea2005}
C.J Read.
\newblock A large weak operator closure for the algebra generated by two
  isometries.
\newblock {\em J. Operator Theory}, 54, no.2:305--316, 2005.

\end{thebibliography}

\end{document}